\documentclass{au}

\usepackage{amssymb,amsfonts,amscd}

\theoremstyle{plain}
\newtheorem{thm}{Theorem}[section]
\newtheorem{lem}[thm]{Lemma}
\newtheorem{prop}[thm]{Proposition}

\theoremstyle{definition}

\newtheorem{defi}[thm]{Definition}

\newtheorem*{rem}{Remark}

\begin{document}

\author{Ewa Graczy\'{n}ska}
\address{Opole University of Technology, Institute of Mathematics
\newline
ul. Luboszycka 3, 45-036 Opole, Poland} \email{egracz@po.opole.pl
\hspace{5mm} http://www.egracz.po.opole.pl/}

\title{On the problem on M-hyperquasivarieties}
\maketitle 

\hspace{6cm} {\it In Memoriam Dietmar Schweigert}

\begin{abstract}
The aim of this paper is to present a solution of the problem 32
posed by K. Denecke and S.L. Wismath in \cite{DW}. It is a
continuation of common results of the author and Dietmar Schweigert
presented in joint papers \cite{gracz22} -- \cite{EGDS5}.

The results of this paper were  partially presented during the conference
Trends in Logic III International in Memoriam Andrzej Mostowski, Helena Rasiowa,
Cecylia Rauszer, Warsaw, September 23, Ruciane-Nida, September 24-25, 2005
and on the Workshop AAA71 and CYA21 in B\c{e}dlewo, Polish Academy of Sciences,
Poland on February 11, 2006.
\end{abstract}


\section{Notation}

Our nomenclature and notation is basically those of G. Birkhoff
\cite{GB1}, K. Denecke and S. Wismath \cite{DW}, G. Gr\"{a}tzer
\cite{GG0},  R. McKenzie, G. McNulty and W. Taylor \cite{ralph}.
Some fundamental concepts and properties of algebras and varieties
may also be found in  \cite{DW}, therefore we omit them here.

\begin{defi}  \label{D.2.1}
{\rm A {\it type}  $\tau$ of an algebra ${\bf A}$ is a function
$\tau : I \rightarrow I\!\!N$ from the indexing set $I$ into the set
$I\!\!N$ of natural numbers, where $\tau (i) = n_{i}$ if $f_{i}$ is
an $n_{i}$-ary operation. A type $\tau$ is finite if the set $I$ is
finite.}
\end{defi}

We deal only with universal algebras of a given $\it{type}$ $\tau :
I \rightarrow I\!\!N$, where $I$ is a nonempty set and $I\!\!N$
denotes the set of all integers.

In that case, for a given $f_{i} \in F$, $\tau(f_{i}) = n_{i}$ is called the
{\it arity} of the operation $f_{i}$, $i \in I$ and we will say that $f_{i}$
is an $n_{i}$-ary operation.

\begin{defi} \label{D.2.2}
{\rm An {\it algebra} \index{algebra} ${\bf A}$ is a pair ${\bf A} = (A,F)$, such that:
$A$ is a nonempty set; $F = (f_{i} : i \in I)$ is a family of
finitary operations in $A$.}
\end{defi}
In the sequel, we shall use the following notation as well:
\begin{defi}  \label{D.2.3}
{\rm An {\it algebra} ${\bf A}$ is a pair ${\bf A} = (A,F^{\bf A})$, such that:
$A$ is a nonempty set; $F = (f_{i}^{\bf A} : i \in I)$ is a set of
finitary operations in $A$.}
\end{defi}

\newtheorem{example}{Example}[section]

\subsection{Identities and varieties of algebras}
By identities of type $\tau$ we mean expressions of  \index{identity}
the form $p \approx q$, where $p,q$ are $n-ary$ polynomial symbols of
a given type $\tau$ for some $n \in I\!\!N$.

A {\it hyperidentity} is formally the same as
identity. An identity $x \approx y$ is called {\it trivial},
where $x$ and $y$ are different variables. \index{trivial identity}

The difference between the concept of an {\it identity} and {\it
hyperidentity} is in {\it satisfaction} (see \cite{gracz22},
\cite{gracz27} and  \cite{DS1}, \cite{DW}).

By $Id(\tau)$ we denote the {\it set of all} identities of type
$\tau$. If $\Sigma$ is a set of identities of type $\tau$, then
$E(\Sigma)$ denotes the closure of $\Sigma$, i.e. the smallest set
of identities of type $\tau$ which contains $\Sigma$ and is closed
under the rules of inference (1)-(5) of G. Birkhoff (see G. Birkhoff
\cite{GB1}, G. Gr\"{a}tzer \cite[p. 170]{GG0},  \cite{ralph} and W.
Taylor  \cite{WT}).

If $V$ is a variety of algebras of type $\tau$ then $Id(V)$ denotes
the set of all {\it identities}, satisfied in $V$, sometimes called
the theory of $V$. If $\Sigma$ is a set of identities of type
$\tau$, then $Mod(\Sigma)$ denotes the {\it class of all models},
i.e. the variety of algebras defined by $\Sigma$.

For a given algebra ${\bf A}$, $Id({\bf A})$ denotes the set of all
identities satisfied in ${\bf A}$. Shortly speaking, a satisfaction
of an identities in an algebra ${\bf A}$ means a satisfaction of a
pair of terms where the variables are bound by universal
quantifiers. A satisfaction of a hyperidentity $p \approx q$ in the
algebra ${\bf A}$ means its satisfaction as the second-order
formula.

{\it Hypersubstitutions} of a given type $\tau$ of terms were
invented by D. Schweigert and the author in \cite{gracz22}. Shortly
speaking, they are mappings sending terms to terms by substituting
variables by (the same) variables and fundamental terms by terms of
the same arities, i.e. $\sigma(x) = x$ for any variable $x$, and for
a given operation symbol $f$, assume that
$\sigma(f(x_{0},...,x_{\tau(f)-1}))$ is a given term of the same
arity as $f$, then $\sigma$ acts on all terms of a given type $\tau$
in an inductive way:
\begin{center}
$\sigma(f(p_{0},...,p_{\tau(f)-1)})) =
\sigma(f(x_{0},...,x_{\tau(f)-1}))(\sigma(p_{0}),...,\sigma(p_{\tau(f)-1}))$,
for $f \in F$,
\end{center}
(where $f(x_{0},...,x_{\tau(f)-1})$ denotes a {\it fundamental
term}).

A different approaches defined the same concept in \cite{lau},
\cite{DW} and \cite{j29}.

$H(\tau)$ denotes the set of all hypersubstitutions $\sigma$ of
a given type $\tau$.

${\bf H}(\tau) = (H(\tau), \circ, \sigma_{id})$ denotes the monoid
of all hypersubstitutions $\sigma$ of a given type $\tau$ with the
operation $\circ$ of composition ant the identity hypersubstitution
$\sigma_{id}$.

\subsection{Quasi-identities and quasivarieties of algebras}

We recall the notion invented by A. I. Mal'cev from \cite{AIM} and \cite{VAG}:
\begin{defi} \label{D.2.6}         \index{quasi-identity}
{\rm A {\it quasi-identity} $e$ is an implication of the form:
\begin{center}
(1.2.1) $(t_{0} \approx s_{0}) \wedge ... \wedge (t_{n-1} \approx
s_{n-1}) \rightarrow (t_{n} \approx s_{n})$.
\end{center}
where $t_{i} \approx s_{i}$ are $k$-ary identities of a given type, for
$i = 0,...,n$.

A {\it quasi-identity} above is {\it satisfied in an algebra} ${\bf
A}$ of a given type if and only if the following implication is
satisfied in ${\bf A}$:  given a sequence $a_{1},...,a_{k}$ of
elements of $A$. If these elements satisfy the equations
$t_{i}(a_{1},...,a_{k}) = s_{i}(a_{1},...,a_{k})$ in ${\bf A}$, for
$i = 0, 1,..., n-1$, then the equality $t_{n}(a_{1},...,a_{k}) =
s_{n}(a_{1},...,a_{k})$ is satisfied in ${\bf A}$. In that case we
write:
\begin{center}
${\bf A} \models (t_{0} \approx s_{0}) \wedge ... \wedge (t_{n-1} \approx
s_{n-1}) \rightarrow (t_{n} \approx s_{n})$.
\end{center}

A quasi-identity $e$ is {\it satisfied in a class} $V$ of algebras
of a given type, if and only if it is satisfied in all algebras
${\bf A}$ belonging to $V$.}
\end{defi}

Following A. I. Mal'cev \cite{AIM} we consider classes $QV$ of
algebras A of a given type $\tau$ defined by quasi-identities and
call them {\it quasivarieties}. \index{quasivariety}  We use the
symbol $QMod(\Sigma)$ for the class $QV$ of algebras of type $\tau$, satisfying
a given set $\Sigma$ of quasi-identities of type $\tau$ and call it the
{\it quasivariety axiomatized} by $\Sigma$. $\Sigma$ is then called a
{\it base} of $QV$.

A {\it hyperquasi-identity} $e$ (of a type $\tau$) is the same as
quasi-identity (of type $\tau$). Sometimes we shall use the notation
for {\it hyperquasi-identity} invented by D. Schweigert in
\cite{DS1}. The difference between quasi-identities and
hyperquasi-identities is in satisfaction. Following the ideas of
\cite{DS1}, part 5, \cite[p. 155]{DW} we modified in \cite{EGDS1}
the definition of the satisfaction of a quasi-identity to the notion
of a hypersatisfaction in the following way:

\begin{defi}  \index{hypersatisfaction}  \label{D.5.3}
{\rm A hyperquasi-identity $e$ is {\it satisfied} (is {\it
hyper-satisfied}, {\it holds}) in an algebra ${\bf A}$ if and only
if the following
implication is satisfied: \\
if $\sigma$ is a hypersubstitution of type $\tau$ and the elements
$a_{1},...,a_{k} \in A$ satisfy the equalities
$\sigma(t_{i})(a_{1},...,a_{k}) = \sigma(s_{i})(a_{1},...,a_{k})$ in
${\bf A}$, for $i \in \{ 0,1,...,k-1 \}$, then the equality
$\sigma(t_{n})(a_{1},...,a_{k}) = \sigma(s_{n})(a_{1},...,a_{k})$
holds in ${\bf A}$.
\begin{center}
In that case, we write $V \models_H e$.
\end{center} }
\end{defi}

In other words, a {\it hyperquasi-identity} is a universally closed
Horn $\forall x \forall \sigma$-formulas, where x varies over all
sequences of individual variables (occurring in terms of the
implication) and $\sigma$ varies over all hypersubstitutions of a
given type. Our modification coincides with Definition 5.1.3 of
\cite{DS1} (see Definition 2.3 of \cite{CCKD}).

\begin{rem}
{\rm All hyperquasi-identities and hyperidentities are written
without quantifiers but they are considered as universally closed
Horn $\forall$-formulas (see \cite{AIM}).}
\end{rem}

\section{Hyperquasi-varieties}

A reformulation of the notion of {\it quasivariety} invented by A. I.
Mal'cev in \cite[p. 210]{AIM} to the notion of {\it hyperquasivariety}
of a given type  $\tau$ was invented by D. Schweigert and the author in
\cite{EGDS1}  in a natural way:
\begin{defi}  \label{D.4.1}
A class $K$ of algebras of type $\tau$ is called a {\it
hyperquasivariety} if there is a set  $\Sigma$ of
hyperquasi-identities of type $\tau$ such that $K$ consists exactly
of those algebras of type $\tau$ that hypersatisfy all the
hyperquasi-identities of $\Sigma$.
\end{defi}

Let us note, that the notion of {\it hyperqasivariety} coincides
with the notion of a {\it hyperquasi-equational class} invented in
\cite[p. 155]{DW}.

\section {Hyperquasi-identities}
We recall only our definitions of \cite{gracz22} of hyperidentities
satisfied in an algebra of a given type and the notion of a {\em
hypervariety}:

\begin{defi}  \label{D.5.1}
An algebra ${\bf A}$ {\it satisfies a hyperidentity} $p \approx q$
if for every hypersubstitution $\sigma \in H(\tau)$ the resulting
identity $\sigma(p) \approx \sigma(q)$ is satisfied in ${\bf A}$. In
this case, we write ${\bf A} \models_{H} p \approx q$. A variety $V$
satisfies a hyperidentity $p \approx q$ if every algebra in the
variety does. In symbols $V \models_{H} p \approx q$.
\end{defi}
\begin{defi}   \label{D.5.2}
{\rm A class $V$ of a algebras of a given type is called a {\it hypervariety}
if and only if $V$ is defined by a set $\Sigma$ of hyperidentities.

In that case we write that $V = HMod(\Sigma)$. }
\end{defi}
\begin{rem}
{\rm Some authors avoid to use our concepts since in \cite{WT4} W. Taylor
defined hypervarieties to be classes of varieties of different types
satisfying certain sets of equations as identities (see also
\cite{NW1}).   }
\end{rem}
The theorem following was proved in \cite{gracz22}:

\begin{thm} \label{T.5.1}
A variety $V$ of type $\tau$ is defined by a set of hyperidentities
if and only if $V = HSPD(V)$, i.e. $V$ is a variety closed under
derived algebras of type $\tau$.
\end{thm}

Let $V$ be a class of algebras of type $\tau$. {\it Derived algebras} were
defined in \cite{PCM}. {\it Derived algebras of a given type} $\tau$ were
defined in \cite{gracz22}.
\begin{defi}  \label{D.5.4}
{\rm Let  ${\bf A} =(A,F)$ be an algebra in $V$ and $\sigma$ a hypersubstitution
in $H(\tau)$. Then the algebra ${\bf B} = (A, (F)^\sigma)$ is a
{\it derived algebra} of ${\bf A}$, with the same universe $A$ and the set
$(F)^\sigma$ of all derived operations of $F$ by $\sigma$.
${\bf B}$ is then denoted as ${\bf A^\sigma}$. }
\end{defi}

$D(V)$ denotes the class of all derived algebras of type $\tau$ of
all algebras of $V$.
\begin{defi}
A quasivariety $V$ is called {\it solid} if and only if $D(V) \subseteq V$.
\end{defi}

In \cite{EGDS1} we presented several theorems of Mal'cev type for solid quasivarieties.

\section{M-solid quasivarieties}

Let $M$ be a subset of $H(\tau)$ closed under compostion $\circ$ and
containing the trivial hypersubstitution i.e. ${\bf M}$ is a
submonoid $(M, \circ, \sigma_{id})$ of the monoid $(H(\tau),\circ,
\sigma_{id})$ .

In \cite{EGDS3} we reformulated the notion of {\em hyperquasivariety} of
\cite{EGDS1} for the case of {\it M-hyperquasivariety} of a given
type in a natural way:
\begin{defi}
{\rm A class $K$ of algebras of type $\tau$ is called an
M-hyperquasi\-variety if there is a set  $\Sigma$ of
M-hyperquasi-identities of type $\tau$ such that $K$ consists
exactly of those algebras of type $\tau$ that M-hyper\-satisfy all
the hyperquasi-identities of $\Sigma$.}
\end{defi}
In \cite[p.155]{DW} M-hyperquasivarieties were called
$M$-hyperquasi-equational classes.

The following follows from \cite{CCKD} (see also \cite{EGDS3}):
\begin{thm} \label{T.11.1}
A quasivariety $K$ of algebras given type is an M-hyper\-quasi\-variety if and
only if it is  M-deriverably closed.
\end{thm}
We accept the following definition of \cite[p. 155]{DW}:
\begin{defi}
{\rm Let $QV$ be a quasivariety, then $QV$ is {\it M-solid} if and only if
every M-derived algebra ${\bf A^\sigma}$ belongs to $QV$, for every
algebra ${\bf A}$ in $QV$ and $\sigma$ in $M$, i.e.
\begin{center}
$D_M(QV) \subset QV$
\end{center} }
\end{defi}

In \cite{EGDS3} we presented some Mal'cev types theorems for {\it
M-hyperquasivarieties}.

\section{M-hyperquasi-identities}

A suitable generalization of our observations made for the set
$H(\tau)$ of all hypersubstitutions was extended in \cite{EGDS3} to
any subset of $H(\tau)$, closed under composition and containing the
trivial hypersubstutuion $\sigma_{id}$. This generalization gives
rise to so called {\em M-hypersubstitutions} of a given type.
\begin{defi}
{\rm $M$-{\it hyperquasi-identity} is formally the same as
quasi-identity.}
\end{defi}
We recall only the definitions of \cite{j29} of the fact that a
hyperidentity is satisfied in an algebra as an M-hyperidentity of a
given type and the notion of M-hypervariety invented in
\cite{gracz27}:

\begin{defi}
{\rm An algebra ${\bf A}$ satisfies a hyperidentity $p \approx q$ as
an M-hyperidentity if for every M-hypersubstitution $\sigma \in M$,
the resulting identity $\sigma(p) \approx \sigma(q)$ holds in ${\bf
A}$.

In that case, we write ${\bf A} \models^{M}_{H} p \approx q$. \\
A variety $V$ satisfies a hyperidentity $p \approx q$ as
M-hyperidentity if every algebra in the variety does. In symbols $V
\models^{M}_{H} p \approx q$. }
\end{defi}
\begin{defi}
{\rm A class $V$ of  algebras of a given type is called an
M-hypervariety if and only if $V$ is defined by a set of
M-hyperidentities.}

In that case we write, that $V = MHMod(\Sigma)$.
\end{defi}
Obviously, the definition above generalizes the notion of a hypervariety to
an M-hypervariety and a hypersatisfaction to an M-hypersatisfaction.
Moreover, every algebra satisfied a set $\Sigma$ as hyperidentities,
satisfies it as a set of of $M$-hyperidentities.
The following was proved in \cite{gracz27}:

\begin{thm}
A variety $V$ of type $\tau$ is defined by a set $\Sigma$ of
M-hyperidentities if and only if $V = HSPD_{M}(V)$, i.e. $V$ is a
variety closed under M-derived algebras of type $\tau$. Moreover, in
this case, the set $\Sigma$ is then M-hypersatisfied in $V$ and $V$
is the class of all M-hypermodels of $\Sigma$, i.e. $V
=MHMod(\Sigma)$.
\end{thm}
In order to explain the difference of the notions invented above with the
notions of Yu. Movsisyan \cite{Mov1}--\cite{Mov4} we invent the following:
\begin{prop}
Let a type $\tau$ and the monoid ${\bf H(\tau)} = (H(\tau), \circ,
\sigma_{id})$ of all hypersubstitutions of type $\tau$ be given.
Then $M_{F}(\tau)$ denotes the set of hypersubstitutions $\sigma$ of
type $\tau$, for which $\sigma(f(x_{0},...,x_{\tau(f)-1}))$ is a
fundamental term and is not a variable, for every functional symbol
$f \in F$. Then the set $M_{F}(\tau)$ is a submonoid of $H(\tau)$.
\end{prop}
\begin{proof} The proof follows from the fact that the composition of two
hypersubstitutions from the set $M_{F}(\tau)$ is a hypersubstiution
in $M_{F}(\tau)$, which is not a projection (i.e. is not a function
determined by a variable).
\end{proof}

\begin{defi}
{\rm Every hypersubstitution of type $\tau$ from the monoid
$M_{F}(\tau)$ is called a $M_{F}(\tau)$-{\it hypersubstitution}. The
monoid ${\bf M_{F}(\tau)} = (M_{F}(\tau), \circ, \sigma_{id})$ is
called the monoid of all $M_{F}(\tau)$-{\it hypersubstitutions} of
type $\tau$.}
\end{defi}
The following shows the connection of the notions invented in
\cite{Mov1}--\cite{Mov4} with the notion of M-hyperidentity of
\cite{gracz22}:
\begin{thm}
For a given algebra ${\bf A}$ of type $\tau$ an identity $p \approx
q$ is satisfied in ${\bf A}$ as a hyperidentity in the sense of
\cite{Mov1}, if and only if it is satisfied in ${\bf A}$ as an
$M_{F}(\tau)$-hyperidentity.
\end{thm}
\begin{defi}
{\rm A hyperquasi-identity $e$  is $M$-{\it hyper-satisfied} ({\it
holds})
in an algebra ${\bf A}$ if and only if the following implication is satisfied:\\
If $\sigma$ is a hypersubstitution of $M$ and the elements
$a_{1},...,a_{n} \in A$ satisfy the equalities
$\sigma(t_{i})(a_{1},...,a_{k}) = \sigma(s_{i})(a_{1},...,a_{k})$ in
${\bf A}$, for $i = 0,1,...,n-1$, then the equality
$\sigma(t_{n})(a_{1},...,a_{k}) = \sigma(s_{n})(a_{1},...,a_{k})$
holds in ${\bf A}$.

We say then, that $e$ is an M-hyperquasi-identity of ${\bf A}$ and
write:
\begin{center}
${\bf A} \models_{H}^{M} (t_{0} \approx s_{0}) \wedge ... \wedge
(t_{n-1} \approx s_{n-1}) \rightarrow (t_{n} \approx s_{n})$.
\end{center}

A hyperquasi-identity $e$  is M-hyper-satisfied (holds) in a class
$V$ if and only if it is M-hypersatisfied in any algebra of $V$. We
write then: $V \models^{M}_{H} e$.  }

\end{defi}
By other words, M-hyperquasi-identity is a universally closed Horn
$\forall x \forall \sigma$-formulas, where x vary over all sequences
of individual variables (occurring in terms of the implication) and
$\sigma$ vary over all hypersubstitutions of $M$. Our modification
coincides with Definition 5.1.3 of \cite{DS1} and Definition 2.3 of
\cite{CCKD}.

\begin{rem}
{\rm All hyperquasi-identities and hyperidentities are written
without quantifiers but they are considered as universally closed
Horn $\forall$-formulas (see \cite{AIM}). In case of a trivial
monoid $M$, the notion of M-hypersatisfaction reduces to the notion
of classical satisfaction of \cite{GB1}, \cite{PCM}. If $M$ is the
monoid of all hypersubstiutions of a given type $\tau$, then the
notion of M-hyperidentity and M-hyperquasi-identity reduces to the
hyperidentity and hyperquasi-identity.}
\end{rem}
\begin{rem}
{\rm Let us note that in case $M$ is a trivial (i.e. 1-element)
monoid of hypersubstitutions of a given type $\tau$, then the
satisfaction $\models^M_H$ gives rise to the satisfaction $\models$
and the operator $D_M$ to the  identity operator.

In case $M = H(\tau)$ we get the notion of $\models_H$ considered in
\cite{EGDS1}.}
\end{rem}

\section{Examples of M-hyperquasi-identities}
\subsection{Quasigroups}
\begin{defi}
{\rm An algebra $(G, \cdot)$ with a binary operation $\cdot$ is called a
{\it quasigroup}, if for all $a \in G$ the operations $x \cdot a$ and
$a \cdot x$ are permutations in $G$.}
\end{defi}
This is equivalent to the fact that in the groupoid ${\bf G}=(G,
\cdot)$ the following two quasi-identities are satisfied:
\begin{center}
(6.1.1) $(x \cdot z \approx y \cdot z) \rightarrow  x\approx y$ and
(6.1.2) $(z \cdot x \approx z \cdot y) \rightarrow  x\approx y$.
\end{center}

\begin{prop}  \label{P.13.1}
If a groupoid ${\bf G}$ satisfies the quasi-identities (6.1.1) and
(6.1.2), then these quasi-identities are satisfied in ${\bf  G}$ as
$M_{3,4}$-hyperquasi-identities, for the monoid ${\bf M_{3,4}} =
(M_{3,4}, \circ, \sigma_{id}) = (\{ \sigma_{3}, \sigma_{4} \},
\circ, \sigma_{id})$, with:
\begin{center}
$M_{3,4} = \{ \sigma_{3}, \sigma_{4} \in H(2): \sigma_{3}(x \cdot y) = x \cdot y, \sigma_{4}(x \cdot y) = y \cdot x \}$.
\end{center}
\end{prop}
\begin{proof}

For $\sigma_{3}(x \cdot y) = x \cdot y$ and $\sigma_{4}(x \cdot y) =
y \cdot x$ the derived quasi-identities:
$\sigma_{3},_{4}(6.1.1,6.1.2)$ are satisfied in ${\bf G}$.
Therefore, the quasi-identities (6.1.1) and (6.1.2) are satisfied in
${\bf G}$ as $M_{3,4}$-hyperquasi-identities in ${\bf G}$.

Note, that the quasi-identities (6.1.1) and (6.1.2) are not
satisfied as $M_{1,2}$-hyperquasi-identities for the monoid
$M_{1,2}$ generated by the first and the second projections:
$\sigma_{1}, \sigma_{2} \in H(2)$. \end{proof}

\subsection{Distributive lattices}
The following proposition is an expression of the example 2 presented in
\cite{TV} in the language of M-hyperidentities:
\begin{prop}
In each distributive lattice ${\bf L} = (L, \wedge, \vee)$ the
following identities holds as $M_{F}(\tau)$-hyperidentities, for the
monoid ${\bf M_{F}(\tau)}$ of all $M_{F}(\tau)$-hypersubstitutions
of distributive lattices:
\begin{center}
{\rm (6.2.1)} $F(x,x) \approx x$;\hspace{5mm}
{\rm (6.2.2)} $F(x,y) \approx F(y,x)$;\\
{\rm (6.2.3)} $F(F(x,y),z) \approx F(x,F(y,z))$; \\
{\rm (6.2.4)} $F(x,G(y,z)) \approx G(F(x,y),G(x,z))$.
\end{center}
\end{prop}
\begin{proof} Let us note, that the monoid $M_{0}(\tau)$ in case of distributive
lattices consists of 4 nonequivalent hypersubstitutions (in the
sense of \cite{JP1}) of type (2,2), namely: $\sigma_{1}(x \wedge y)
= x \wedge y$, $\sigma_{1}(x \vee y) = x \vee y$, $\sigma_{2}(x
\wedge y) = x \vee y$, $\sigma_{2}(x \vee y) = x \wedge y$,
$\sigma_{3}(x \wedge y) = x \wedge y = \sigma_{3}(x \vee y)$,
$\sigma_{4}(x \vee y) = x \vee  y= \sigma_{4}(x \vee y)$. The
hypersubstitutions of all the identities of (6.2.1) -- (6.2.2) are
identities satisfied in any distributive lattice ${\bf L}$.
\end{proof}

\subsection{Boolean Algebras}

We express example 1 of \cite{TV} in the language of M-hyperidentities:
\begin{prop}
The following identity holds as an $M_{F}(\tau)$-hyperiden\-tity in
every Boolean algebra:
\begin{center}
$F(G(x,y)',z)' \approx G(F(x',z)',F(y',z)')$.
\end{center}
\end{prop}

\subsection{Flat algebras}

{\it Flat algebras} were invented by R. McKenzie. They were considered in
\cite{JMK} as specific 0-{\it smilattice algebras}.

Let $\sigma$ be a finite   signature containing (among other symbols)
a binary symbol $\wedge$ (the meet) and a nullary symnol 0.
\begin{defi}
{\rm By a 0-{\it semilattice  $\tau$-algebra} we mean an algebra of type $\tau$
satisfying the equations

(6.5.1) $x \wedge (y \wedge z) \approx (x \wedge y) \wedge z$;

(6.5.2) $x \wedge y \approx y \wedge x$;

(6.5.3) $x \wedge x \approx x$;

(6.5.4) $f(x_{1},...,x_{i-1},0,x_{i+1},...,x_{n}) \approx 0$,

for every n-ary operation $f$ of type $\tau$ and every $i \in \{1,...,n \}$.}
\end{defi}
\begin{defi}
{\rm A {\it flat algebra} is a 0-semilattice algebra ${\bf A}$ such that
$a \wedge b = 0$ for all pairs of distinct elements $a,b \in A$.}
\end{defi}
Consider the monoid ${\bf M}_{0,\wedge}(\tau)$ of all prehypersubstitutions
of type $\tau$ leaving the constant 0 and the operation $\wedge$ unchanged.

Then the following holds:
\begin{thm}
The variety of flat algebras is ${\bf M}_{0,\wedge}(\tau)$-solid.
\end{thm}
\begin{proof} Given a 0-semilattice (flat algebra) ${\bf A}$ and a
hypersubstitution $\sigma \in M(0,\wedge)$. Then obviously the
derived identities of identities (6.5.1) and (6.5.3) remains
unchanged and satisfied in ${\bf A}$. Consider the derived identity
of (6.5.4) by $\sigma$, i.e.
$\sigma(f(x_{1},...,x_{i-1},0,x_{i+1},...,x_{n}) ) \approx \sigma(0)
\approx 0$, i.e. sis satified in ${\bf A}$. Consider the derived
identity of (6.5.4) by $\sigma$, i.e.
$\sigma(f)(x_{1},...,x_{i-1},0,x_{i+1},...,x_{n}) ) \approx 0$,
which is satisfied in ${\bf A}$ as $\sigma(f)$ is an n-ary
polynomial symbol $g$ of type $\tau$. Moreover, in every
0-semilattice (flat algebra) the following equation holds:
\begin{center}
(6.5.4*) $p(x_{1},...,x_{i-1},0,x_{i+1},...,x_{n}) \approx 0$,
\end{center}
for every nontrivial (i.e. not a variable) term $p$. We prove this
fact by induction on the complexity of the term $p$, which is not a
variable. Assume that the induction hypothesis holds for n-ary
polynomials $p_{1},...,p_{m}$ and let $p = g(p_{1},...,p_{m})$ for
an m-ary functional symbol $g$. \\
Then $g(p_{1},...,p_{m})(x_{1},...,x_{i-1},0,x_{i+1},...,x_{n})
\approx g(0,...,0) \approx 0$. \end{proof}

\begin{defi}  \index{compatible 0-semilattice}
{\rm A 0-semilattice $\tau$-algebra is {\it compatible} if it satisfies the
equation:

(6.5.5) $f(z_{1},...,z_{i-1},x \wedge y,z_{i+1},...,z_{n}) \approx$

          $f(z_{1},...,z_{i-1},x,z_{i+1},...,z_{n}) \wedge
          f(z_{1},...,z_{i-1},y,z_{i+1},...,z_{n})$,

for every n-ary operation $f$ of type $\tau$ and every $i \in \{1,...,n \}$.}
\end{defi}
Consider the monoid ${\bf M}^{*}_{0,\wedge}(\tau)$ of all prehypersubstitutions
of type $\tau$ leaving the constant 0 and the operation $\wedge$ unchenged in
such a way, that $\sigma(f)$ is always a functional symbol (of the same arity
as $f$), for every $f$ of type $\tau$. Then the following is obvious:
\begin{thm}
The variety of 0-semilattice algebras is ${\bf M}^{*}_{0,\wedge}(\tau)$-solid.
The variety of compatible flat algebras is ${\bf M}^{*}_{0,\wedge}(\tau)$ solid.
\end{thm}
Recall from \cite[p. 666]{JMK} the following definition of {\it basic x-term}
of depth $n$:     \index{basic x-term}
\begin{defi}
{\rm The term $x$ is the only {\it basic x-term} of depth 0. For $n
> 0$, {\it basic x-terms} of depth $n$ are the terms
$f(x_{1},...,x_{i-1},t,x_{i+1},...,x_{n})$ such that $f$ is an
$n$-ary operation symbol of type $\tau$, $1 \leq i \leq n$, $t$ is a
basic $x$-term of depth $n-1$ and $x_{1},...,x_{n}$ are variables
different of $x$.}
\end{defi}
\begin{lem} \label{L.13.1}
For every hypersubstitution $\sigma$ from ${\bf M}^{*}_{0,\wedge}(\tau)$
and every basic $x$-term $t(x)$ of depth $n$, the hypersubstitution term
$\sigma(t(x))$ is a basic $x$-term of depht $n$.
\end{lem}
\begin{proof} The term $\sigma(f(x_{1},...x_{i-1},t,x_{i+1},...,x_{n}))$
equals to a term of the form:\\ $g(x_{1},...,x_{i-1},\sigma(t),x_{i+1},...,x_{n})$
for some $n$-ary functional symbol $g$ ot type $\tau$. Therefore the lemma follows
by induction on the complexity of a basic $x$-term $t(x)$. \end{proof}

Recall from \cite[p. 668]{JMK}, that for a finite compatible, flat
algebra ${\bf A}$ there exists a finitely q-based (i.e. having a
finite base of quasi-identities) quasivariety $Q'_{A}$ containing
\index{q-based} ${\bf A}$. The base constructed contains all
identities of the form (6.5.1) -- (6.5.5) and quasi-identities
constracted by means of basic $x$- and $y$-terms of depth $\leq K$,
for $K$ being the cardinality of $A$ and the operation $\wedge$. Via
Lemma  \ref{L.13.1} we conclude the following slight strengthening
of Lemma 3.1 of \cite{JMK}:
\begin{prop}
$Q'_{A}$  is a finitely q-based ${\bf M}_{0,\wedge}(\tau)$-hyperquasivariety
containing ${\bf A}$.
\end{prop}

\section{Hyperquasi-equational logic}

In this section we present a solution of the following particular
case of the Problem 32 \cite[p. 291]{DW}:

(P.32) {\it Give the derivation rules for M-hyperquasi-equational
logic}.

First, we shall consider the case where the monoid $M$ is trivial,
i.e. one-element. \index{equational logic} In the sequel, ${\bf E}$
denotes the {\it equational logic}, i.e. the fragment of the
first-order logic without relation symbols. The formulas of ${\bf
E}$ are all possible identities of a given type $\tau$, the set of
axioms $Eq$ of ${\bf E}$ are identities of the form $p \approx p$,
and the rules of inferences are the equality rules (atomic formulas
are regarded as identities) and the {\it substitution rule}, i.e. G.
Birkhoff's rules (1)--(5) of derivation. \\$E$ denotes the set of
equality axioms of a given type $\tau$ (see \cite[p. 33]{VAG}).

For a set $\Sigma$ of (hyper)quasi-identities of a given type
$\tau$, $HQMod(\Sigma)$ denotes the class of all algebras ${\bf A}$
which hypersatisfy all elements of $\Sigma$. \index{hyperequational
logic}

${\bf HE}$ denotes the {\it hyperequational logic}, i.e. the
fragment of the second-order logic, extending the equational logic.
The formulas and axioms are the same as in ${\bf E}$. To the
inference rules we add the rule (6) of hypersubstitution defined in
\cite{gracz22}.

Following \cite[p. 72]{VAG} a quasi-identity $e$ is called a {\it
consequence} of the set $\Sigma$ of quasi-identities if for every
algebra ${\bf A}$ of type $\tau$, ${\bf A} \models \Sigma$ implies
that ${\bf A} \models e$. In symbols: $\Sigma \models e$.

We say that an identity $e$ is a {\it hyperconsequence} of a set of
quasi-identities $\Sigma$, if for every algebra ${\bf A} \in
HMod(\Sigma)$, it follows that ${\bf A} \models_{H} e$, i.e. ${\bf
A} \models _{H} \Sigma$ implies ${\bf A} \models_{H} e$. In symbols:
$\Sigma \models_{H} e$.

Following \cite{VAG} we use the following notation:
\begin{center}
$\Delta \rightarrow \alpha$, for a set $\Delta = \{ p_{i} \approx
q_{i} : 0 \leq i \leq n-1 \}$ and $\alpha = p_{n} \approx q_{n}$
\end{center}
instead of the quasi-identity:
\begin{center}
$p_{0} \approx q_{0} \wedge ... \wedge p_{n-1} \approx q_{n-1}
\rightarrow p_{n} \approx q_{n}$.
\end{center}

We adopt the convention, that an identity $p \approx q$ may be
regarded as a quasi-identity $e$ of the form $\emptyset \rightarrow
p \approx q$, where $\emptyset$ denotes the empty set.

G. Birkhoff's well known theorem is called {\it the completeness
theorem}:
\begin{thm}   \index{completeness theorem}
An identity $e$ is a consequence of a set $\Sigma$ of identities if
and only if $e$ is derived from $\Sigma$ in ${\bf E}$.
\end{thm}
The question naturally arises of when an identity is a consequence
of a set of quasiidentities $\Sigma$ (see \cite{GB1}). Following
\cite[p. 72]{VAG} it is necessary, together with a substitution rule
to consider the {\it modus ponens} rule:  \index{modus ponens}
\begin{center}
(MP)  $\frac{\alpha, \{ \alpha \} \cup \Delta  \rightarrow
\beta}{\Delta \rightarrow \beta}$.
\end{center}
\index{quasi-equational logic}

Recall from \cite[p. 73]{VAG}, that in the quasi-equational logic
${\bf Q}$ (of a given type $\tau$), without relation symbols, the
formulas are all possible quasi-identities of a given type $\tau$,
the axioms are the {\it equality axioms} (E.1) -- (E.4) and the
inference rules are the {\it substitution rule}, the {\it cut rule}
and the {\it extension rule}. We list all of them.

Axioms:\\ (E.1) the reflexivity:
\begin{center}
$p \approx p \rightarrow p \approx p$,
\end{center}
(E.2) the symmetry: \index{symmetry axiom}
\begin{center}
$p \approx q \rightarrow q \approx p$,
\end{center}
(E.3) the transitivity: \index{transitivity axiom }
\begin{center}
$(p \approx q) \wedge (q \approx r) \rightarrow (p \approx r)$,
\end{center}
or in an  equivalent notation:
\begin{center}
$\{ p \approx q, q \approx r \} \rightarrow p \approx r$,
\end{center}
(E.4) the compatibility: \index{compatibility axiom}
\begin{center} $(t_{0} \approx s_{0})
\wedge ... \wedge (t_{\tau(f)-1} \approx s_{\tau(f)-1}) \rightarrow
(f(t_{0},...,t_{\tau(f)-1}) \approx f(s_{0},...,s_{\tau(f)-1}))$,
\end{center}
for every operation symbol $f$ of type $\tau$,

or in an equivalent notation:
\begin{center}
$\{ t_{0} \approx s_{0},...,t_{\tau(f)-1} \approx s_{\tau(f)-1} \}
\rightarrow (f(t_{0},...,t_{\tau(f)-1}) \approx
 f(s_{0},...,s_{\tau(f)-1}))$,
\end{center}
for every operation symbol $f$ of type $\tau$.

The inference rules are the following rules: \index{substitution
rule}

(7.1) the substitution rule (where $\delta$ is a substitution of
variables):
\begin{center}
$\frac{ \{ \gamma_{0},...,\gamma_{n-1} \} \rightarrow \beta}{ \{
\delta(\gamma_{0}),...,\delta(\gamma_{n-1}) \} \rightarrow
\delta(\beta)}$
\end{center}

(7.2) the cut rule:      \index{cut rule}
\begin{center}
$\frac{\Delta \rightarrow \alpha, \{ \alpha \} \cup \Gamma
\rightarrow \beta} {\Delta \cup \Gamma \rightarrow \beta}$
\end{center}

(7.3) the extension rule:      \index{extension rule}
\begin{center}
$\frac{\Delta \rightarrow \alpha}{ \{ \beta \} \cup \Delta
\rightarrow \alpha}$.
\end{center}

We write $\Sigma \vdash_{Q} e$ if there exists a derivation of a
quasi-identity $e$ from a set $\Sigma$ of quasi-identities in ${\bf
Q}$.

The classical result by many authors is the following:
\begin{thm}
A quasi-identity $e$ is a consequence of a set $\Sigma$ of
quasi-identities  if and only if $e$ is derivable from $\Sigma$ in
${\bf Q}$.
\begin{center}
In symbols: $\Sigma \models_{Q} e$ if and only if $\Sigma \vdash_{Q}
e$.
\end{center}
\end{thm}
We modify {\em quasi equational logic} {\bf Q} to {\it
hyperquasi-equational logic} {\bf HQ} by adding a new rule:

\index{hypersubstitution rule} (7.4) a hypersubstitution rule (where
$\sigma$ is a hypersubstitution of $H(\tau)$):
\begin{center}
$\frac{(t_{0} \approx s_{0}) \wedge ... \wedge (t_{n-1} \approx
s_{n-1}) \rightarrow (t_{n} \approx s_{n})}{\sigma(t_{0}) \approx
\sigma(s_{0}) \wedge ... \wedge \sigma(t_{n-1}) \approx
\sigma(s_{n-1}) \rightarrow \sigma(t_{n}) \approx \sigma(s_{n})}$,
\end{center}

or in an equivalent notation:

(7.4) a hypersubstitution rule (where $\sigma$ is a
hypersubstitution of $H(\tau)$):
\begin{center}
$\frac{ \{ \gamma_{0},...,\gamma_{n-1} \} \rightarrow \beta}{ \{
\sigma(\gamma_{0}),...,\sigma(\gamma_{n-1}) \} \rightarrow
\sigma(\beta)}$.
\end{center}

\begin{defi} \label{D.7.1} \index{hyperquasi-equational logic}
{\rm By {\bf HQ} we denote the hyperquasi-equational logic, which is
an extension of the hyperequational logic {\bf HE} based on the
equality axioms $E$ and four rules (7.1) -- (7.4) above.}
\end{defi}

We write $\Sigma \vdash_{HQ} e$ if there exists a derivation of $e$
from $\Sigma$ in ${\bf HQ}$. \\
We write $\Sigma \models_{HQ} e$ if $e$ is a hyperconsequence of
$\Sigma$, considered as a hyperbase, i.e. if ${\bf A} \in
HQMod(\Sigma)$, then ${\bf A} \models_{HQ} e$.
\begin{defi}
{\rm A set $\Sigma$ of quasi-identities of type $\tau$ is called
{\it hyperclosed} if and only if it is closed under the equality
axioms, the substitution rule, hypersubstitution rule, the cut rule,
the extension rule. } \index{hyperclosed set of quasi-identities}
\end{defi}
We reformulate the classical results in the following way:
\begin{thm} \label{T.7.1}
A set $\Sigma$ is a set of all (hyper)quasi-identities of type
$\tau$, (hyper)satisfied in a class $K$ of algebras of type $\tau$
if and only if it is (hyper)closed.
\end{thm}
{\it Proof}. If $\Sigma$ is a set of all hyperquasi-identities
hypersatisfied in a class $K$ of algebras of type $\tau$, then it is
closed in ${\bf Q}$, i.e. is closed under the equality axioms and
the substitution rule, the cut and the extension rule. In
consequence it is also closed under the rules of equational logic.
If $e$ is a quasi-identity of $\Sigma$, then for every $\sigma \in
H(\tau)$, the hypersubstitution $\sigma(e)$ of $e$ by $\sigma$ is
satisfied in $K$. Therefore $\Sigma$ is closed under the
hypersubstitution rule (7.4). In case if $e$ is an identity of type
$\tau$, we conclude that $\sigma(e)$ is satisfied in $K$ for every
$\sigma \in H(\tau)$. Therefore $\Sigma$ is closed under the rule
(6) of hypersubstitution, i.e. is hyperclosed.

Assume now, that $\Sigma$ is hyperclosed. Therefore it is closed. We
conclude that $\Sigma$ is a set of quasi-identities satisfied in a
class $K$ of algebras of type $\tau$. As $\Sigma$ is hyperclosed,
therefore for every quasi-identity $e$ of $\Sigma$ and every $\sigma
\in H(\tau)$, the derived quasi-identity $\sigma(e)$ is also
satisfied by $K$, which means that $K$ is a class of algebras of
type $\tau$, which hypersatisfies $\Sigma$. $\Box$

The clue of the next proofs is the following:

\begin{prop}    \label{P.7.1}
A derivation from $\Sigma$ in {\bf HQ} means a derivation from
$(7.4)(\Sigma)$ in {\bf Q}, i.e. one first need to close the set
$\Sigma$ under the hypersubstitution rule {\rm (7.4)} and then under
the equality axioms and other rules. The resulting set will be
already closed under all axioms and inference rules of ${\bf HQ}$.
\end{prop}
More precisely:
\begin{prop} \label{P.7.com}
The hypersubstitution rule {\rm (7.4)} commutes with all the axioms
and rules of the logic {\bf HQ}.
\end{prop}
{\it Proof}. First, we note that the assertion easily holds the
equality axioms (E1)--(E3). Moreover, by an easy induction on the
complexity of terms, the following generalization of the rule (E.4)
is valid in the logic ${\bf Q}$:
\begin{center}
(GE.4) $\{ t_{0} \approx s_{0},...,t_{\tau(f)-1} \approx
s_{\tau(f)-1} \} \rightarrow(p(t_{0},...,t_{\tau(f)-1}) \approx
p(s_{0},...,s_{\tau(f)-1}))$,
\end{center}
for every term $p$ of type $\tau$.

We prove that if the axiom (E.4) is applied first:
\begin{center}
$ \{ t_{0} \approx s_{0},...,t_{\tau(f)-1} \approx s_{\tau(f)-1} \}
\rightarrow (f(t_{0},...,t_{\tau(f)-1}) \approx
f(s_{0},...,s_{\tau(f)-1}))$,
\end{center}
and then the hypersubstitution rule (7.4) is applied by a
hypersubstitution $\sigma$:

\begin{center}
$ \{ \sigma(t_{0}) \approx \sigma(s_{0}),...,\sigma(t_{\tau(f)-1})
\approx \sigma(s_{\tau(f)-1}) \} \rightarrow
(\sigma(f(t_{0},...,t_{\tau(f)-1})) \approx
\sigma(f(s_{0}),...,(s_{\tau(f)-1})))$,
\end{center}
then one may apply rule (GE.4) with $p =
\sigma(f(x_{0},...,x_{n}))$, to obtain the resulting quasi-identity:
$ \{ \sigma(t_{0}) \approx \sigma(s_{0}),...,\sigma(t_{\tau(f)-1})
\approx \sigma(s_{\tau(f)-1}) \} \rightarrow$\\$ \rightarrow
(\sigma(f)(\sigma(t_{0}),...,\sigma(t_{\tau(f)-1})) \approx
\sigma(f)(\sigma(s_{0}),...,\sigma(s_{\tau(f)-1})))$.

Now we prove the assertion for the modus ponens rule (MP):
\begin{center}
(MP)  $\frac{\alpha, \{ \alpha \} \cup \Delta \rightarrow
\beta}{\Delta \rightarrow \beta}$.
\end{center}
i.e. we will show, that if the (MP) rule is applied first and then
the hypersubstitution rule (7.4) is applied to deduce a
quasi-identity $e = \sigma (\Delta) \rightarrow \sigma (\beta)$,
then one may apply the hypersubstitution rule (7.4) first to
$\alpha$ and $\alpha \cup \Delta \rightarrow \beta$ and then (MP),
which leads to the quasi-identity $e$ as well.

Secondly, assume that the substitution rule (7.1) is applied (where
$\delta$ is a substitution of variables):
\begin{center}
(7.1) $\frac{ \{ \gamma_{0},...,\gamma_{n-1} \} \rightarrow \beta}{
\{ \delta(\gamma_{0}),...,\delta(\gamma_{n-1}) \} \rightarrow
\delta(\beta)}$
\end{center}
and then the hypersubstitution rule (7.4) is applied to get the
quasi-identity:
\begin{center}
(*) $\{ \sigma (\delta(\gamma_{0})),...,\sigma (\delta
(\gamma_{n-1})) \} \rightarrow \sigma (\delta (\beta))$
\end{center}
for some hypersubstitution $\sigma \in H(\tau)$ and a substitution
$\delta$ of variables. Assume that the substitution $\delta$ acts on
variables $x_{0},...,x_{m}$ of terms
$\gamma_{0},...,\gamma_{n-1},\beta$ putting: $\delta(x_{k})=p_{k}$,
then putting $\delta_{1}(x_{k})=\sigma(p_{k})$ on variables of terms
$\sigma(p_{k})$ of type $\tau$, we get that:
$\sigma(\delta_{1}(\gamma_{i})) = \delta_{1}(\sigma(\gamma_{i}))$,
for $i = 0,...,n-1$
and $\sigma(\delta_{1}(\beta)) = \delta_{1}(\sigma(\beta))$. \\
We conclude that the quasi-identity (*) equals to the
quasi-identity:
\begin{center}
(*) $\{ \delta_{1}(\sigma (\gamma_{0})),...,\delta_{1}(\sigma
(\gamma_{n-1})) \} \rightarrow \delta_{1}(\sigma (\beta))$,
\end{center}
which means that one may apply the hypersubstitution rule (7.4)
first and then the substitution rule (7.1) to get the same result.

The proof for the cut rule is similar. Assume that the cut rule
(7.2) is applied:

\begin{center}
(7.2) $\frac{\Delta \rightarrow \alpha, \{ \alpha \} \cup \Gamma
\rightarrow \beta} {\Delta \cup \Gamma \rightarrow \beta}$
\end{center}
and then the hypersubstitution rule (7.4) by a hypersubstitution
$\sigma$ gives rise to the quasi-identity:
\begin{center}
(**) $\sigma (\Delta) \cup \sigma (\Gamma) \rightarrow \sigma
(\beta)$.
\end{center}
Then one may apply the hypersubstitution rule (7.4) by $\sigma$ to
the quasi-identities:
\begin{center}
$\Delta \rightarrow \alpha$ and $ \{ \alpha \} \cup \Gamma
\rightarrow \beta$
\end{center}
to get the resulting quasi-identity  (**) via the cut rule (7.2).

We finalize with the proof of the statement for the extension rule,
applying first:
\begin{center}
(7.3) $\frac{\Delta \rightarrow \alpha}{ \{ \beta \} \cup \Delta
\rightarrow \alpha}$
\end{center}
and assuming that the hypersubstitution rule (7.3) by $\sigma$ was
applied then, leading to the quasi-identity:
\begin{center}
(***) $\{ \sigma (\beta) \} \cup \sigma (\Delta) \rightarrow \sigma
(\alpha)$.
\end{center}
Then one may apply the hypersubstitution rule (7.4) $\sigma$ first
to the quasi-identity: $\Delta \rightarrow \alpha$, to get the
resulting quasi-identity (***) as a result of the extension rule
(7.3). $\Box$

The observation above is a generalization of that which has been
already noticed in \cite[p. 121]{gracz17}, for the fact that
derivation rules (1)-(5) of G. Birkhoff and the new rule (6) of
hypersubstitution behave similarly, i.e. closing a set $\Sigma$ of
identities under (1)-(6) means, to close $\Sigma$ under (6) first
and then under rules (1)-(5) and we are done.

Therefore, we can say that the hyperquasi-equational logic is the
one-step extension of the quasi-equational logic by the
hypersubstitution rule (7.4).

We obtain a slight generalization of Corollary 2.2.3 of \cite[p.
72]{VAG}:
\begin{prop} \label{P.7.2}
An identity $e$ is a (hyper)consequence of a set $\Sigma$ of
quasi-identities if and only if there is a derivation of $e$ (of
$\sigma(e)$, for every $\sigma \in H(\tau)$)  from $E \cup \Sigma$
by the substitution rule and modus ponens rule (and the
hypersubstitution rule {\rm (7.4)}).
\end{prop}

{\it Proof}. The first part of the theorem for ${\bf Q}$ is the
classical result (see \cite{VAG}).

Assume that an identity $e$ is a hyperconsequence of a set $\Sigma$,
i.e. $\Sigma \models_{HQ} e$. It means, that for every algebra ${\bf
A}$ if ${\bf A} \models_{HQ} \Sigma$, then ${\bf A} \models_{HQ} e$.
In other words: for every algebra ${\bf A}$ if ${\bf A} \models_{Q}
(7.4)(\Sigma) = \{ \sigma(\Sigma): \sigma \in H(\tau) \}$, then
${\bf A} \models_{Q} \sigma(e)$, for every $\sigma \in H(\tau)$.
Therefore, we conclude that $(7.4)(\Sigma) \models_{Q} \sigma(e)$,
for every $\sigma \in H(\tau)$. Therefore, via Corollary 2.2.3 of
\cite[p. 72]{VAG}, we conclude, that for every $\sigma \in H(\tau)$
there is a derivation of $\sigma(e)$ from $E \cup \Sigma$ by the
substitution rule and the modus ponens rule.

Assume now, that there is a derivation of $e$ from $E \cup \Sigma$
by the substitution, hypersubstitution and modus ponens rule. Then
for every $\sigma \in H(\tau)$ there is a derivation of $\sigma(e)$
from $E \cup \Sigma$ by the substitution, hypersubstitution and
modus ponens rule. Applying the proposition \ref{P.7.1}, we conclude
that there is a derivation of $\sigma(e)$  from the closure $(7.4)(E
\cup \Sigma)$ of the set $E \cup \Sigma$ by (7.4), by the
substitution and modus ponens rule, for every $\sigma \in H(\tau)$.
By Corollary 2.2.3 of \cite[p. 72]{VAG}, we conclude  that
$\sigma(e)$ is a consequence of $(7.4)\Sigma$, for every $\sigma \in
H(\tau)$, i.e. $(7.4)\Sigma \models_{Q} \sigma(e)$, for every
$\sigma \in H(\tau)$. Therefore $\Sigma \models_{HQ} e$. $\Box$

The following is the modification of the classical {\it completeness
theorem} of the logic {\bf Q}:
\begin{thm} \label{T.7.2}
A (hyper)quasi-identity $e$ is a (hyper)consequence of a set
$\Sigma$ of (hyper)quasi-identities if and only if it is derivable
from $\Sigma$ in {\bf (H)Q}.
\begin{center}
In symbols: $\Sigma \models_{(H)Q} e$ if and only if $\Sigma
\vdash_{(H)Q} e$.
\end{center}
\end{thm}

{\it Proof}. The part of the theorem for {\bf Q} is the classical
result of Selman \cite{SA} (see \cite[p. 73]{VAG}).

Assume that $\Sigma \models_{HQ} e$, i.e. if an algebra ${\bf A} \in
HQMod(\Sigma)$, i.e. if ${\bf A} \models_{HQ} \Sigma$, then ${\bf A}
\models_{HQ} e$. This is equivalent to the implication: if ${\bf A}
\models_{Q} \sigma(\Sigma)$, for every $\sigma \in H(\tau)$, then
${\bf A} \models_{HQ} e$. Equivalently we write this implication in
the following way: if ${\bf A} \models_{Q} (7.4)(\Sigma)$, then
${\bf A} \models_{HQ} e$. From the completeness theorem of G.
Birkhoff theorem 2.2.5 \cite[p. 73]{VAG} for the logic ${\bf Q}$, we
conclude, that if $(7.4)(\Sigma) \models_{Q} \sigma(e)$, for every
$\sigma \in H(\tau)$, i.e. if ${\bf A} \models_{Q} (7.4)(\Sigma)$,
then ${\bf A} \models_{Q} \sigma(e)$,
for every $\sigma \in H(\tau)$. Therefore we conclude the implication:\\
$\Sigma \models_{HQ} e$, then $(7.4)(\Sigma) \vdash \sigma(e)$, for
every $\sigma \in H(\tau)$.  We got: $\Sigma \vdash_{HQ} e$.

Assume now that $e$ is derivable  from $\Sigma$ in ${\bf HQ}$, i.e.
$\Sigma \vdash_{HQ} e$. By proposition \ref{P.7.1} we conclude that
the quasi-identity $e$ is derivable from $(7.4)(\Sigma)$ in ${\bf
Q}$, i.e. $(7.4)(\Sigma) \vdash_{Q} e$. Therefore, via completeness
theorem for ${\bf Q}$, we obtain that $(7.4)(\Sigma) \models_{Q} e$,
i.e. for every algebra ${\bf A}$, such that ${\bf A} \models
(7.4)(\Sigma)$ it follows that ${\bf A} \models e$. This means, that
from ${\bf A} \models_{HQ} \Sigma$ it follows that ${\bf A} \models
e$. The similar argument follows for every derived quasi-identity
$\sigma(e)$, of $e$, for every $\sigma \in H(\tau)$. Namely, if
$\Sigma \vdash_{HQ} e$, then for every $\sigma \in H(\tau)$ we
conclude, that $\Sigma_{HQ} \sigma(e)$, as if $e_{1},...,e_{n}$ is a
proof of $e$ from $\Sigma$ in {\bf HQ}, then:
$e_{1},\sigma(e_{1}),...,\sigma(e_{n})$ is a proof of $\sigma(e)$
from $\Sigma$ in {\bf HQ}. Therefore we got $\Sigma \vdash_{HQ}
\sigma(e)$, for every $\sigma \in H(\tau)$.  Finally we conclude
that $\Sigma \models_{HQ} e$.  $\Box$

\section{M-hyperquasi-equational logic}

In this section we present a solution  of the Problem 32 \cite[p. 291]{DW}:

(P.32) {\it Give the derivation rules for M-hyperquasi-equational logic}.

Assume that a monoid ${\bf M} = (M, \circ, \sigma_{id})$ of
hypersubstitutions of type $\tau$ is given.

For a set $\Sigma$ of (hyper)quasi-identities of a given type
$\tau$, $MHQMod(\Sigma)$ denotes the class of all algebras ${\bf A}$
which hypersatisfy all elements of $\Sigma$.
\index{M-hyperequational logic} By ${\bf MHE}$ we denote the M-{\it
hyperequational logic}, i.e. the fragment of the second-order logic,
without relation symbols, extending the equational logic. The
formulas and axioms are the same as in ${\bf E}$. To the inference
rules of ${\bf E}$ we add the rule $(6)_{M}$ of M-hypersubstitution
defined by the author in \cite{gracz27}.

We modify {\em quasi equational logic} ${\bf Q}$ to M-{\it
hyperquasi-equational logic} {\bf MHQ} by adding a new rule:
\index{M-hypersubstitution rule}

(8.4) an M-{\it hypersubstitution rule} (where $\sigma$ is a
hypersubstitution of M):
\begin{center}
$\frac{(t_{0} \approx s_{0}) \wedge ... \wedge (t_{n-1} \approx
s_{n-1}) \rightarrow (t_{n} \approx s_{n})} {\sigma(t_{0}) \approx
\sigma(s_{0}) \wedge ... \wedge \sigma(t_{n-1}) \approx
\sigma(s_{n-1}) \rightarrow \sigma(t_{n}) \approx \sigma(s_{n})}$,
\end{center}

or in an equivalent notation:

(8.4) an M-{\it hypersubstitution rule} (where $\sigma$ is a
hypersubstitution of M):
\begin{center}
$\frac{ \{ \gamma_{0},...,\gamma_{n-1} \} \rightarrow \beta}{ \{
\sigma(\gamma_{0}),...,\sigma(\gamma_{n-1}) \} \rightarrow
\sigma(\beta)}$
\end{center}

The rule above generalizes the hypersubstitution rule (7.4) of {\bf
HQ}.

\begin{defi}
{\rm By {\bf MHQ} we denote the M-hyperquasi-equational logic, which
is an extension of M-hyperequational logic {\bf MHE}, generalizes
the logic ${\bf HQ}$ and is based on the equality axioms
(E.1)--(E.4) and the three inference rules of (7.1) - (7.3) of {\bf
Q} and the M-hypersubstitution rule (8.4).}
\end{defi}

We write $\Sigma \vdash^{M}_{HQ} e$ if there exists a derivation of
$e$ from $\Sigma$ in ${\bf MHQ}$.

We write $\Sigma \models^{M}_{HQ} e$ if $e$ is an M-hyperconsequence
of $\Sigma$, considered as a hyperbase, i.e. if ${\bf A} \in
MHQMod(\Sigma)$, then ${\bf A} \models^{M}_{HQ} e$.

\begin{rem}
{\rm Note, that if $M$ is a trivial monoid, then the logic ${\bf
MHQ}$ coincides with the logic ${\bf Q}$.  If $M = H(\tau)$, then
the logic {\bf MHQ} coincides with the logic {\bf HQ}. }
\end{rem}
\begin{defi}
{\rm A set of quasi-identities of type $\tau$ is called M-{\it
hyperclosed} if and only if it is closed under the equality axioms
and the substitution rule, M-hypersubsti\-tution rule, the cut rule
and the extension rule.}
\end{defi}
We generalize the classical results in the following way:
\begin{thm} \label{T.17.1}
A set $\Sigma$ is a set of all (M-hyper)quasi-identities of a class
$K$ of algebras of type $\tau$ if and only if it is M-hyperclosed.
\end{thm}

{\it Proof}. If $\Sigma$ is a set of (M-hyper)quasi-identities
M-hypersatisfied in a class $K$ of algebras of type $\tau$, then it
is closed in ${\bf Q}$, i.e. is closed under the equality axioms and
the substitution rule, the cut and the extension rule. In
consequence it closed under the rules of equational logic, i.e. G.
Birkhoff's rules (1) - (5). If $e$ is a quasi-identity of $\Sigma$,
then for every $\sigma \in M$, the hypersubstitution $\sigma(e)$ of
$e$ by $\sigma$ is satisfied. Therefore $\Sigma$ is closed under the
M-hypersubstitution rule (8.4). In case if $e$ is an identity, we
conclude that $\sigma(e)$ is satisfied for every $\sigma \in M$.
Therefore $\Sigma$ is closed under the rule $(6)_{M}$ of
M-hypersubstitution (which is a particular case of the rule (8.4)),
i.e. under {\bf MHE}.

Assume now, that $\Sigma$ is M-hyperclosed. Therefore it is closed.
We conclude, that $\Sigma$ is a set of quasi-identities satisfied in
a class $K$ of algebras of type $\tau$. As $\Sigma$ is
M-hyperclosed, therefore for every quasi-identity $e$ of $\Sigma$
and every $\sigma \in M$, the derived quasi-identity $\sigma(e)$ is
also satisfied by $K$, which means that $K$ is a class of algebras
of type $\tau$, which M-hypersatisfy $\Sigma$. $\Box$

The clue of the next proofs is the following common generalization
of Proposition \ref{P.7.1}:

A derivation from $\Sigma$ in {\bf MHQ} means a derivation from
$(8.4)(\Sigma)$ in {\bf Q}, i.e. one first need to close the set
$\Sigma$ under the hypersubstitution rule (8.4) and then under the
equality axioms and the other rules. The resulting set will be
already closed under all inference rules of ${\bf MHQ}$.

More precisely:
\begin{prop} \label{P.17.1}
The M-hypersubstitution rule {\rm (8.4)} commutes with all the
axioms and rules of the logic {\bf MHQ}.

In symbols: for a quasi-identity $e$, the following equivalence
holds:
\begin{center}
$\Sigma \vdash ^{M}_{HQ} e$ if and only if $(8.4)(\Sigma) \vdash
_{Q} e$.
\end{center}
\end{prop}

{\it Proof}. By Proposition \ref{P.7.com} the assertion holds for
the set $E$ of equality axioms, as in the rule (8.4) one should
consider hypersubstitutions $\sigma \in M$ only. We prove the
assertion for the modus ponens rule (MP):
\begin{center}
(MP)  $\frac{\alpha, \{ \alpha \cup \Delta \} \rightarrow
\beta}{\Delta \rightarrow \beta}$.
\end{center}
i.e. we will show, that if the (MP) rule is applied first and then
the M-hypersubstitution rule (8.4) is applied to deduce a
quasi-identity $e = \sigma (\Delta) \rightarrow \sigma (\beta)$,
then one may apply the M-hypersubstitution rule (8.4) first to
$\alpha$ and $\alpha \cup \Delta \rightarrow \beta$ and then (MP),
which leads to the quasi-identity $e$ as well.

Secondly, assume that the substitution rule (7.1) is applied (where
$\delta$ is a substitution of variables):
\begin{center}
(7.1) $\frac{ \{ \gamma_{0},...,\gamma_{n-1} \} \rightarrow \beta}{
\{ \delta(\gamma_{0}),...,\delta(\gamma_{n-1}) \} \rightarrow
\delta(\beta)}$
\end{center}
and then the hypersubstitution rule (8.4) is applied to get the
quasi-identity:
\begin{center}
(*) $\{ \sigma (\delta(\gamma_{0})),...,\sigma (\delta
(\gamma_{n-1})) \} \rightarrow \sigma (\delta (\beta))$
\end{center}
for some M-hypersubstitution $\sigma \in M$ and a substitution
$\delta$ of variables. Asume that the substitution $\delta$ acts on
variables $x_{0},...,x_{m}$ of terms
$\gamma_{0},...,\gamma_{n-1},\beta$ putting: $\delta(x_{k})=p_{k}$,
then putting $\delta_{1}(x_{k})=\sigma(p_{k})$ on variables of terms
$\sigma(p_{k})$ of type $\tau$, we get that:
$\sigma(\delta_{1}(\gamma_{i})) = \delta_{1}(\sigma(\gamma_{i}))$,
for $i = 0,...,n-1$
and $\sigma(\delta_{1}(\beta)) = \delta_{1}(\sigma(\beta))$. \\
We conclude that the quasi-identity (*) is equal to the
quasi-identity:
\begin{center}
(*) $\{ \delta_{1}(\sigma (\gamma_{0})),...,\delta_{1}(\sigma
(\gamma_{n-1})) \} \rightarrow \delta_{1}(\sigma (\beta))$,
\end{center}
which means that one may apply the M-hypersubstitution rule (8.4)
first and then the substitution rule (7.1) to get the same result.

The proof for the cut rule is similar. Assume that the cut rule
(7.2) is applied:

\begin{center}
(7.2) $\frac{\Delta \rightarrow \alpha, \{ \alpha \} \cup \Gamma
\rightarrow \beta} {\Delta \cup \Gamma \rightarrow \beta}$
\end{center}
and then the M-hypersubstitution rule (8.4) by a hypersubstitution
$\sigma \in M$ gives rise to the quasi-identity:
\begin{center}
(**) $\sigma (\Delta) \cup \sigma (\Gamma) \rightarrow \sigma
(\beta)$.
\end{center}
Then one may apply the hypersubstitution rule (8.4) by $\sigma$ to
the quasi-identities:
\begin{center}
$\Delta \rightarrow \alpha$ and $ \{ \alpha \} \cup \Gamma
\rightarrow \beta$
\end{center}
to get the resulting quasi-identity  (**) via the cut rule (7.2).

We finalize with the proof of the statement for the extension rule,
applied first:
\begin{center}
(7.3) $\frac{\Delta \rightarrow \alpha}{ \{ \beta \} \cup \Delta
\rightarrow \alpha}$
\end{center}
and assuming that the M-hypersubstitution rule (8.4) by $\sigma$ was
applied then, leading to the quasi-identity:
\begin{center}
(***) $\{ \sigma (\beta) \} \cup \sigma (\Delta) \rightarrow \sigma
(\alpha)$.
\end{center}
Then one may apply the M-hypersubstitution rule (8.4) $\sigma$ first
to the quasi-identity: $\Delta \rightarrow \alpha$, to get the
resulting quasi-identity (***) as a result of the extension rule
(7.3). $\Box$

The observation above is a generalization of that the author has
already noticed in \cite{gracz17}, for the fact that derivation
rules (1)-(5) of G. Birkhoff and the new rule $(6)_{M}$ of
hypersubstitution behave similarly, i.e. closing a set $\Sigma$ of
identities under $(1)-(6)_{M}$ means, to close $\Sigma$ under
$(6)_{M}$ first and then under rules (1)-(5) and we get that the
resulting set is closed under the rules $(1)-(6)_{M}$.

Therefore, we can say that the M-hyperquasi-equational logic ${\bf
MHQ}$ is the one-step extension of the quasi-equational logic {\bf
Q} by the M-hyper\-substi\-tution rule (8.4).

We obtain a slight generalization of Corollary 2.2.3 of
\cite[p.72]{VAG}:
\begin{prop} \label{P.17.2}
An identity e is an (M-hyper)consequence of a set of
quasi-identities $\Sigma$ if and only if there is a derivation of
$e$ ($\sigma(e)$, for every $\sigma \in M$), from $E \cup \Sigma$ by
the substitution rule and modus ponens rule (and the
M-hypersubstitution rule $(8.4)$).
\end{prop}
{\it Proof}. Assume that an identity $e$ is an M-hyperconsequence of
a set $\Sigma$, i.e. $\Sigma \models^{M}_{HQ} e$. It means, that for
every algebra ${\bf A}$ if ${\bf A} \models^{M}_{HQ} \Sigma$, then
${\bf A} \models^{M}_{HQ} e$. In other words: for every algebra
${\bf A}$ if ${\bf A} \models_{Q} (8.4)(\Sigma)$, then ${\bf A}
\models_{Q} \sigma(e)$, for every $\sigma \in M$. Therefore, we
conclude that $(8.4)(\Sigma) \models_{Q} \sigma(e)$, for every
$\sigma \in M$. Therefore, via Corollary 2.2.3 of \cite[p. 72]{VAG},
we conclude, that for every $\sigma \in M$ there is a derivation of
$\sigma(e)$ from $(8.4)(E \cup \Sigma)$ by the substitution rule and
the modus ponens rule, i.e. from the set $E \cup \Sigma$ by the
substitution rule, the modus ponens rule and the M-hypersubstitution
rule.

Assume now, that there is a derivation of $e$ from $E \cup \Sigma$
by  substitution, M-hypersubstitution and modus ponens rule. Then
for every $\sigma \in M$ there is a derivation of $\sigma(e)$ from
$E \cup \Sigma$ by the substitution, M-hypersubstitution and modus
ponens rule. Applying Proposition \ref{P.17.1}, we conclude that
there is a derivation of $\sigma(e)$  from the closure $(8.4)(E \cup
\Sigma)$ of the set $E \cup \Sigma$ by (8.4), by the substitution
and modus ponens rule, for every $\sigma \in M$. By Corollary 2.2.3
of \cite[p. 72]{VAG}, we conclude  that $\sigma(e)$ is a consequence
of $(8.4)\Sigma$, for every $\sigma \in M$, i.e. $(8.4)\Sigma
\models_{Q} \sigma(e)$, for every $\sigma \in M$. Therefore $\Sigma
\models^{M}_{HQ} e$. $\Box$
\begin{rem}
{\rm Note, that in case of a trivial $M$, Proposition \ref{P.17.2}
is nothing else but Corollary 2.2.3 of \cite[p. 72]{VAG}. }
\end{rem}
\begin{thm} \label{T.17.2}
A (hyper)quasi-identity $e$ is an M-hyperconsequence of a set
$\Sigma$ of (hyper)quasi-identities if and only if it is derivable
from $\Sigma$ in {\bf MHQ}.
\begin{center}
In symbols: $\Sigma \models^{M}_{HQ} e$ if and only if $\Sigma
\vdash^{M}_{HQ} e$.
\end{center}
\end{thm}
{\it Proof}. Assume that $\Sigma \models^{M}_{HQ} e$, i.e. if an
algebra ${\bf A} \in MHQMod(\Sigma)$, i.e. if ${\bf A} \models
^{M}_{HQ} \Sigma$, then ${\bf A} \models^{M}_{HQ} e$. This is
equivalent to the implication: if ${\bf A} \models_{Q}
\sigma(\Sigma)$, for every $\sigma \in M$, then ${\bf A}
\models^{M}_{HQ} e$. Equivalently we write this implication in the
following way: if ${\bf A} \models_{Q} (8.4)(\Sigma)$, then ${\bf A}
\models^{M}_{HQ} e$. From the completeness theorem of G. Birkhoff
Theorem 2.2.5 \cite[p. 73]{VAG} for the logic ${\bf Q}$, we
conclude, that if $(8.4)(\Sigma) \models_{Q} \sigma(e)$, for every
$\sigma \in M$, i.e. if ${\bf A} \models_{Q} (8.4)(\Sigma)$, then
${\bf A} \models_{Q} \sigma(e)$, for every $\sigma \in M$, i.e.
$(8.4)(\Sigma) \models _{Q} \sigma(e)$, for
every $\sigma \in M$.\\
Therefore we conclude the implication:\\
$\Sigma \models^{M}_{HQ} e$, then $(8.4)(\Sigma) \vdash \sigma(e)$,
for every $\sigma \in M$.  We got: $\Sigma \vdash^{M}_{HQ} e$.

Assume now that $e$ is derivable  from $\Sigma$ in ${\bf MHQ}$, i.e.
$\Sigma \vdash^{M}_{HQ} e$. By proposition \ref{P.17.1} we conclude
that the quasi-identity $e$ is derivable from $(8.4)(\Sigma)$ in
${\bf Q}$, i.e. $(8.4)(\Sigma) \vdash_{Q} e$. Therefore, via
completeness theorem for ${\bf Q}$, we obtain that $(8.4)(\Sigma)
\models_{Q} e$, i.e. for every algebra ${\bf A}$, such that ${\bf A}
\models (8.4)(\Sigma)$ it follows that ${\bf A} \models e$. This
means, that from ${\bf A} \models^{M}_{HQ} \Sigma$ it follows that
${\bf A} \models e$. The similar argument follows for every derived
quasi-identity $\sigma(e)$, for every $\sigma \in M$. Namely, if
$\Sigma \vdash ^{M}_{HQ} e$, then for every $\sigma \in M$ we
conclude, that $\Sigma \vdash^{M}_{HQ} \sigma(e)$, as if
$e_{1},...,e_{n}$ is a proof of $e$ from $\Sigma$ in {\bf MHQ},
then: $e_{1},\sigma(e_{1}),...,\sigma(e_{n})$ is a proof of
$\sigma(e)$ from $\Sigma$ in {\bf MHQ}. Therefore we got, that
$\Sigma \vdash \sigma(e)$. Finally we conclude that $\Sigma \models
^{M}_{HQ} e$. $\Box$


\begin{thebibliography}{100}




\bibitem{AJ}
Acz\'{e}l, J., {\em Proof of a theorem of distributive type hyperidentities},
Algebra Universalis {\bf 1}, 1971, 1--6.

\bibitem{TCB}
Bartee, T. C., {\em Digital Computer Fundamentals}, McGraw-Hill, 1966.

\bibitem{GB0}
Bartee, T. C., Birkhoff, G., {\em Modern Applied Algebra}, Corrected
third printing, Mc-Graw Hill Book Company, 1970.

\bibitem{BV}
Belousov, V. D., {\em Systems of quasigroups with generalized
identities}, Uspechi Mat. Nauk {\bf 20}, 1965, 75--146. English
translation: Russian Math. Surveys 20, 75--143.

\bibitem{WB1}
Belousov, W.,
{\em Introduction to the theory of quasigroups and loops} (in Russian),
Moscow, 1967.

\bibitem{bergman}
Bergaman, G. M., {\em  Hyperidentities of grous and semigroups}, Aequationes
Mathematicae {\bf 23}, 1981, 50--65.

\bibitem{GB1}
Birkhoff, G.,
{\em  On the structure of abstract algebras}, Proc. Cambr.
Philos. Soc. {\bf 31}, 1935, 433--454.

\bibitem{burris}
Burris, S., Sankappanavar, H. P.,
{\em A Course in Universal Algebra}, Springer Verlag, New York, 1981.

\bibitem{CCKD}
Chompoonut, Ch., Denecke, K. {\em M-solid Quasivarieties}, East-West
J. of Mathematics, Vol. {\bf 4}, No. 2, 2002, 177--190.

\bibitem{PCM}
Cohn, P.M., \emph{Universal Algebra}, Reidel, 1981, Dordreht.

\bibitem{DW}
Denecke, K.,  Wismath, S.L., {\em Hyperidentities and clones}, Algebra, Logic
and Applications Vol. 14, Gordon and Breach Science Publishers, 2000.
ISBN90-5699-235-X. ISSN: 1041-5394.

\bibitem{lau}
Denecke, K., Lau, D.,  P\"{o}schel, R., Schweigert, D., {\em
Hyperidentities, hyperequational classes and clone congruences}, in:
Contributions to General Algebra {\bf 7},
Verlag-H\"{o}lder-Pichler-Tempsky, Wien, 1991, 97--117.

\bibitem{VAG}
Gorbunov, V.A., {\em Algebraic Theory of Quasivarieties},
Consultants Buereau, 1998, New York, USA.

\bibitem{gracz17}
Graczy\'{n}ska, E.,
{\em On normal and regular identities and hyperidentities}, Proceedings of the
V Universal Algebra Symposium, Universal and Applied Algebra, Turawa, Poland,
3--7 May, 1988, World Scientific, 1989, 107--135.

\bibitem{gracz27}
Graczy\'{n}ska, E., {\em G. Birkfoff's theorems for $M$-solid
varieties}, Algebra Universalis {\bf 40}, 1998, pp. 109--117.

\bibitem{EG1}
Graczy\'{n}ska, E., {\em On the problem of basis for
hyperquasivarieties}, in: Contributions to General Algebra 16,
Proceedings of the Dresden Conference 2004  (AAA 68) and the Summer
School 2004, Verlag Johannes Heyn, Klagenfurt 2005, 91--98.

\bibitem{gracz22}
Graczy\'{n}ska, E., Schweigert, D.,
{\em Hyperidentities of a given type}, Algebra Universalis {\bf 27}, 1990,
305--318.

\bibitem{EGDS1}
Graczy\'{n}ska, E.,  Schweigert, D., {\em Hyperquasivarieties},
Preprint Nr. 336, ISSN 0943-8874, August 2003, Fachbereich Mathematik,
Universit\"{a}t Kaiserslautern (Germany).

\bibitem{EGDS2}
Graczy\'{n}ska, E., Schweigert, D. {\em Hybrid bases for varieties of
semigrous}, Algebra Universalis {\bf 50}, 2003, 129--139.


\bibitem{EGDS3}
Graczy\'{n}ska, E. Schweigert, D., {\em M-hyperquasivarieties}, Demonstratio
Mathematica, Vol. XXXIX, No. 1, 2006, 33--42.

\bibitem{EGDS4}
Graczy\'{n}ska, E. Schweigert, D., {\em Derived and fluid varieties},
in print.

\bibitem{EGDS5}
Graczy\'{n}ska, E., Schweigert, D., {\em The dimension of a
variety}, in print in Discussiones Mat.

\bibitem{GG0}
Gr\"{a}tzer, G.,
{\em Universal Algebra}, 1st ed., D. Van Nostrand Company, Inc.,
printed in the USA, 1968.

\bibitem{GG1}
Gr\"{a}tzer, G., {\em Universal Algebra}, 2nd ed., Springer Verlag,
Berlin, 1979.

\bibitem{DHRK}
Hobby, D.,  McKenzie, R., \emph{The Structure of finite algebras},
AMS, Vol. 76, Providence Rhode Island, USA, 1988.

\bibitem{JMK}
Je\u{z}ek, J., Mar\'{o}ti, M., McKenzie, R., \emph{ Quasiequational theories
of flat algebras}, Czechoslovak Mathematical Journal {\bf 55} (130), 2005,
665-675.

\bibitem{ralph}
McKenzie, R.,  McNulty, G., Taylor, W., {\em Algebras, Lattices, Varieties},
Vol. {\bf 1}, 1987.

\bibitem{KDJK}
Koppitz, J., Denecke, K., \emph{M-solid varieties of algebras},
Advances in Mathematics, Vol. {\bf 10}, Springer, 2006.


\bibitem{AIM1}
Mal'cev, A. I.,
{\em Multiplication of classes of algebraic systems} (in Russian),
Sibirskij Math. J. {\bf 8}, 1967, 346--365.

\bibitem{AIM}
Mal'cev, A. I., \emph{Algebraic systems}, Springer Verlag Berlin
Heidelberg New York, 1973.

\bibitem{MM}
Mar\'{o}ti, M. (see Je\u{z}ek, J., Mar\'{o}ti, M., McKenzie, R. \cite{JMK})

\bibitem{Mov1}
Movsisyan, Yu. M., {\em Introduction to the theory of algebras with
hyperidentities} (Russian), Izdat. Erevan. Univ., Erevan, 1986.

\bibitem{Mov2}
Movsisyan, Yu. M., {\em Hyperidentities and hypervarieties in Algebras}
(in Russian), Izdat. Erevan. Univ., Erevan, 1990.

\bibitem{Mov3}
Movsisyan, Yu. M., {\em Hyperidentities of Boolean algebras} (in Russian)
Izv. Ross. Akad. Nauk, Ser. Math. {\bf 56}, 1992, no. 3, 654--672.

\bibitem{Mov4}
Movsisyan, Yu. M., {\em Algebras with hyperidentities of variety of Boolean algebras} (in Russian)
Izv. Ross. Akad. Nauk, Ser. Math. {\bf 60}, 1996, No. 6, 127--168.

\bibitem{Mov5}
Movsisyan, Yu. M. {\em Hypersidentities and hypervarieties}, Scientiae
Mathematicae  Japonicae, {\bf 54}, 2001, 595--640.  ISSN 1346-0447.

\bibitem{neumann}
Neumann, W. D.,
{\em On Malcev conditions}, J. Austral. Math. Soc. {\bf 1}, 1974, 376--384.

\bibitem{NW1}
Neumann, W. D., {\em Representing varieties of algebras by algebras},
J. Austral. Math. Soc. {\bf 11}, 1970, 1--8.

\bibitem{NW}
Neumann, W. D., {\em Mal'cev conditions, spectrs and Kronecker product},
J. Austral. Math. Soc. (A), {\bf 25}, 1987, 103--117.

\bibitem{nulty}
McNulty, G.,
(see McKenzie, R., McNulty G.,  Taylor, W., \cite{ralph})

\bibitem{PP1}
Penner, P.,
{\em  Hyperidentities of semilattices}, Houston J. of Math. {\bf 10}, 1984,
81--108.

\bibitem{JP1}
P{\l}onka, J., \emph{Proper and inner hypersubstitutions of varieties},
Proceedings of the International Conference Summer School on General
Algebra and Ordered Sets, Olomouc 1994, eds. I. Chajda, R. Hala\u{s},
F. Krutsk\'{y}, 106--116.

\bibitem{j29}
P{\l}onka, J.,
{\em On hyperidentities of some varieties}
in: General Algebra and Descrete Mathematics, eds: K. Denecke, O. L\"{u}ders,
Heldermann Verlag Berlin, 1995, 199-213.

\bibitem{BMS1}
Schein, B. M.,
{\em On the theory of generalized groups} (in Russian),
Dokl. Acad. Nauk SSSR 153, 1963, 296--299. (M.R. 30, No. 1200).

\bibitem{DS1}
Schweigert, D., \emph{Hyperidentities}, in: I. G. Rosenberg and G.
Sabidussi, Algebras and Orders, Kluwer Academic Publishers, 1993,
405--506. ISBN 0-7923-2143-X.

\bibitem{DS2}
Schweigert, D., \emph{On derived varieties}, Discussiones Mathematicae
Algebra and Stochastic Methods {\bf 18}, 1998, 17--26.

\bibitem{SA}
Selman, A., \emph{Completeness  of calculi for axiomatically defined classes
of algebras}, Algebra Universalis {\bf 2}, No. 1, 20--32, 1972.

\bibitem{tarski}
Tarski, A.,
{\em A remark on functionally free algebras}, Ann. of Math. {\bf 47}, 1946,
163--167.

\bibitem{WT}
Taylor, W.,
{\em Equational logic}, Houston J. Math. {\bf 5}, 1979, 1--83.

\bibitem{WT3}
Taylor, W.,
{\em Characterizing Mal'tsev conditions}, Algebra Universalis {\bf 3}, 1973,
351--384.

\bibitem{WT4}
Taylor, W.,
{\em Hyperidentities and hypervarieties}, Aequationes Math. {\bf 23}, 1981,
30--49.

\bibitem{WT2}
(see McKenzie, R., McNulty, G., Taylor, W., \cite{ralph})

\bibitem{TV}
Toghanyan, W., {\em Subdirectly irreducible algebras with various
equations}, Aequationes Math. {\bf 68},  2004, 98--107.

\bibitem{SW1}
Wismath, S.L., {\em Unary hyperidentities for type $<1>$ algebras},
Discussiones Math. Vol. {\bf 17}, No. 1, 1997, 105--112.

\bibitem{SLW}
Wismath, S.L.,
(see Denecke, K. and Wismath, S.L. \cite{DW})
\end{thebibliography}
\end{document}